\documentclass{article}
\usepackage{amsmath,amsfonts,amssymb,amsthm,latexsym,stmaryrd}

\allowdisplaybreaks[4]

\newcommand{\st}{\mathrel{|}}

\newcommand{\K}{\mathcal{K}}		

\newcommand{\FFF}{\mathcal{F}}		

\DeclareMathOperator{\III}{\mathbb{I}}		

\newcommand{\bbbp}{\mathbb{P}}			
\DeclareMathOperator{\Prob}{\bbbp}
\DeclareMathOperator{\UpperProb}{\lefteqn{\smash{\overline{\bbbp}}}\phantom{\bbbp}}
\DeclareMathOperator{\LowerProb}{\underline{\bbbp}}	

\newcommand{\bbbr}{\mathbb{R}}		
\newtheorem{lemma}{Lemma}
\newtheorem{proposition}{Proposition}
\newtheorem{corollary}{Corollary}
\newtheorem{theorem}{Theorem}

\newlength{\IndentI}
\newlength{\IndentII}
\newlength{\IndentIII}
\newlength{\IndentIV}
\newlength{\WidthI}
\newlength{\WidthII}
\newlength{\WidthIII}
\newlength{\WidthIV}
\title{The generality of the zero-one laws}
\author{Akimichi Takemura$^1$\\
  {\small takemura@stat.t.u-tokyo.ac.jp, http://www.e.u-tokyo.ac.jp/\~{}takemura}\\[1mm]
  Vladimir Vovk$^2$\\
  {\small vovk@cs.rhul.ac.uk, http://vovk.net}\\[1mm]
  Glenn Shafer$^{3,2}$\\
  {\small gshafer@andromeda.rutgers.edu, http://glennshafer.com}}

\theoremstyle{definition}

\date{August 2009}
\begin{document}
\maketitle
\begin{abstract}
  We prove game-theoretic generalizations of some well known zero-one laws.
  Our proofs make the martingales behind the laws explicit, and our results illustrate 
  how martingale arguments can have implications going beyond measure-theoretic probability.
\end{abstract}


\footnotetext[1]{Department of Mathematical Informatics,
  Graduate School of Information Science and Technology,
  University of Tokyo,
  7-3-1 Hongo, Bunkyo-ku, Tokyo 113-0033, Japan.}
\footnotetext[2]{Computer Learning Research Centre,
  Department of Computer Science,
  Royal Holloway, University of London,
  Egham, Surrey TW20 0EX, England.}
\footnotetext[3]{Department of Accounting and Information Systems,
  Rutgers Business School---Newark and New Bruns\-wick,
  180 University Avenue, Newark, NJ 07102, USA.}

\section{Introduction}
\label{sec:introduction}

Using simple martingale arguments, we generalize three zero-one laws:
\begin{description}
\item[Kolmogorov's zero-one law:]
  In an infinite sequence of independent trials,
  an event whose happening or failing
  is not affected by a finite number of the trials
  has either probability zero or probability one.
  This was first proven by Kolmogorov in an appendix
  to his \emph{Grundbegriffe}, published in 1933. 
\item[Ergodicity of Bernoulli shifts:]
  In an infinite sequence of independent and identically distributed trials,
  an event that is invariant under
  shifts has either probability zero or probability one.
  This remains true for an infinite sequence of states of a homogeneous Markov chain,
  under certain conditions.
  See, for example, 
  Cornfeld et al.\ (1982), \S 8.1, Theorems 1 and 2.
\item[Hewitt and Savage's zero-one law:]
  In an infinite sequence of independent and identically distributed trials,
  a permutable event has either probability zero or probability one.
  This was first proven by Hewitt and Savage, in 1955
  (Hewitt and Savage (1955), Theorem 11.3).
\end{description}
As one would expect from previous literature
(cf.\ the bibliographical notes to Chapter 8 of 
Cornfeld et al.\ (1982)),
our argument for the ergodicity of Bernoulli and Markov shifts
is similar to but simpler than our argument
for Kolmogorov's zero-one law.
In the case of Hewitt and Savage's zero-one law,
our game-theoretic generalization concerns only a special case;
it is an open question whether it holds for all permutable events.

We formulate our generalizations in the game-theoretic framework
introduced by Shafer and Vovk (2001)
\nocite{takeuchi:2004}.
This means that we use martingale ideas only, without resorting to measure-theoretic 
assumptions.  Our arguments thus make clear that the zero-one laws have a 
martingale meaning that extends beyond the measure-theoretic setting.
Instead of beginning 
with a probability measure that determines prices (expected values)
for all measurable and bounded payoffs,
the game-theoretic framework begins directly with prices.
As soon as prices are given, we have martingales:
a martingale is the capital process resulting from a strategy
for gambling at the given prices.
As soon as we have martingales, we can prove theorems by constructing strategies.
For example, we can prove that an event has probability one
by constructing a strategy that multiplies the capital it risks
by an infinite factor if the event fails.

If enough prices are given to determine a probability measure,
then each event will have a probability,
and the capital process for a gambling strategy will be a martingale
in the measure-theoretic sense.
But we do not assume that the prices determine a probability
measure, and so in general events have only upper and lower probabilities.  Thus
our concept of a martingale is more general than the measure-theoretic concept, and
our results are more general than the classical zero-one laws.

As we will see, even a limited number of prices
will determine for each subset $E$ of the sample space
a lower probability $\LowerProb(E)$
and an upper probability $\UpperProb(E)$ satisfying
$
  0
  \le
  \LowerProb(E)
  \le
  \UpperProb(E)
  \le
  1
$
and
$
  \LowerProb(E)
  =
  1
  -
  \UpperProb
  \left(
    E^c
  \right),
$
where $E^c$ is the complement of $E$.
(The latter equality will be our definition of lower probability.)
We consider $E$ certain if $\LowerProb(E)=1$
or, equivalently, $\UpperProb(E^c)=0$.
We consider $E$ impossible if $E^c$ is certain.
In the special case where there are enough prices
to determine a probability measure on a $\sigma$-algebra containing $E$,
$\LowerProb(E)$ and $\UpperProb(E)$
will both equal the probability the measure gives to $E$.

In the usual theory,
a zero-one law specifies a property of an event $E$
that guarantees that either $\Prob(E)=0$ or $\Prob(E)=1$.
The corresponding game-theoretic zero-one law might say that
either $\LowerProb(E)=0$ or $\UpperProb(E)=1$.
This does not assert that anything is certain or impossible,
but it reduces to the usual law
in the case where the prices determine a probability measure
on a $\sigma$-algebra containing $E$,
for then upper and lower probabilities are probabilities.
Sometimes we will be able to prove a stronger statement,
such as the disjunction of $\UpperProb(E)=0$, $\LowerProb(E)=1$,
and $0=\LowerProb(E)<\UpperProb(E)=1$.

We will first lay out protocols for our games (\S\ref{sec:gtp}),
define upper and lower probabilities in these protocols (\S\ref{sec:upplowprob}),
and explain what a zero-one law looks like
in terms of upper and lower probabilities (\S\ref{sec:twotypes}).
Then we give our martingale proofs:
first for Kolmogorov's zero-one law for tail events
in independently priced trials (\S\ref{sec:kolmogorov}),
for ergodicity in independently and identically priced trials
and Markov trials (\S\ref{sec:ergodic}),
and finally for Hewitt and Savage's zero-one law for permutable events
in independently and identically priced trials (\S\ref{sec:hewitt-savage}).

\section{Protocols}
\label{sec:gtp}

In the protocols considered in this article,
two players, whom we call Skeptic and Reality,
play an infinite number of rounds, which we call \emph{trials}.
On each trial,
Skeptic chooses a gamble and then Reality determines its payoff.
Each player sees the other's moves as they are made.
Reality sees Skeptic's move before determining its payoff,
and Skeptic sees Reality's move before they go on to the next trial.

Formally, Skeptic chooses a real-valued function $F$ on a set $\Omega$,
Reality then chooses an element $\omega$ of $\Omega$,
and Skeptic's payoff for the trial is $F(\omega)$.
We call $\omega$ the \emph{outcome} of the trial.
We write $\K_n$ for Skeptic's capital after the $n$th trial,
and we assume that his initial capital is one monetary unit ($\K_0 = 1$).

We assume that Skeptic chooses $F$ from a non-empty set $\FFF$
of real-valued functions on $\Omega$ with these three properties:
\begin{enumerate}
\item
  If $F_1\in\FFF$ and $F_2\in\FFF$, then $F_1+F_2\in\FFF$.
\item
  If $c\ge0$ and $F\in\FFF$, then $cF\in\FFF$.
\item
  There is no $F\in\FFF$ such that $F(\omega)>0$ for all $\omega\in\Omega$.
\end{enumerate}
Property~2 guarantees that the function on $\Omega$ that is identically equal
to $0$ is in $\FFF$.
Properties~1 and~2 are the defining properties of a \emph{cone}; we call property~3
\emph{coherence}.

We illustrate properties 1--3 in the simple case of $\Omega=\{0,1\}$.
Skeptic chooses a pair of real numbers
$(F(0),F(1))\in\FFF \subseteq \mathbb{R}^2$,
where $\FFF$ is a convex cone in $\mathbb{R}^2$
(properties 1 and 2).
By the coherence property
$\FFF$ cannot meet the strictly positive quadrant $(0,\infty)^2$.
Conversely, any convex cone $\FFF$ not meeting $(0,\infty)^2$
satisfies properties 1--3.

If $P$ is a probability measure on a $\sigma$-algebra for $\Omega$,
then the set of all measurable real-valued functions on $\Omega$
that have expected value zero with respect to $P$ is a coherent cone.
Let us call it the \emph{zero cone} for $P$. The set of all 
measurable real-valued functions on $\Omega$
that have nonpositive expected value with respect to $P$ is also a coherent cone;
we may call it the \emph{nonpositive cone} for $P$.
But a coherent cone need not be a zero cone or a nonpositive cone.
Coherent cones that are proper subsets of zero cones are 
studied in Shafer and Vovk (2001). 
A cone of this type may include a function
$F$ such that $-F$ is not in the cone; this means that Skeptic can take only 
one side of the gamble represented by $F$.
Because our results apply to any coherent cone, they 
generalize the measure-theoretic zero-one laws.

We consider three different protocols,
which differ only in how $\FFF$ may change from one trial to the next.
In the first protocol, $\FFF$ is always the same.
If $\FFF$ is the zero cone for a probability measure $P$,
this is a protocol for betting on successive outcomes at odds given by $P$.
If $\FFF$ is the nonpositive cone for a probability measure $P$,
it is a protocol for betting on successive outcomes at odds 
no better than those given by $P$.

\medskip

\noindent
\textsc{Protocol 1. Identically priced trials}

\noindent
\textbf{Parameters:} set $\Omega$; coherent cone $\FFF$ of real-valued functions on $\Omega$

\noindent
\textbf{Protocol:}

\parshape=6
\IndentI   \WidthI
\IndentI   \WidthI
\IndentII  \WidthII
\IndentII  \WidthII
\IndentII  \WidthII
\IndentI   \WidthI
\noindent
$\K_0 := 1$.\\
FOR $n=1,2,\ldots$:\\
  Skeptic announces $F_n\in\FFF$.\\
  Reality announces $\omega_n\in\Omega$.\\
  $\K_n := \K_{n-1} + F_n(\omega_n)$.\\
END FOR

\medskip

In the second protocol,
the cone from which Skeptic selects his move may change from trial to trial,
but the cone $\FFF_n$ for the $n$th trial is fixed at the beginning of the game;
it does not depend on outcomes of previous trials.  In the special case where 
each $\FFF_n$ is the zero cone (resp., nonpositive cone)
for a probability measure, this is a protocol for betting on successive trials
at odds that are fair (resp., not favorable) according
to probability measures assigned to the trials at the outset of the game.

\medskip

\noindent
\textsc{Protocol 2. Independently priced trials}

\noindent
\textbf{Parameters:} set $\Omega$;
  coherent cones $\FFF_1,\FFF_2,\ldots$ of real-valued functions on $\Omega$

\noindent
\textbf{Protocol:}

\parshape=6
\IndentI   \WidthI
\IndentI   \WidthI
\IndentII  \WidthII
\IndentII  \WidthII
\IndentII  \WidthII
\IndentI   \WidthI
\noindent
$\K_0 := 1$.\\
FOR $n=1,2,\ldots$:\\
  Skeptic announces $F_n\in\FFF_n$.\\
  Reality announces $\omega_n\in\Omega$.\\
  $\K_n := \K_{n-1} + F_n(\omega_n)$.\\
END FOR

\medskip

In the third protocol, the cone for the $n$th trial may depend
on the outcome $\omega_{n-1}$ of the preceding trial.
To indicate this possible dependence,
we designate the cone by $\FFF(\omega_{n-1})$.
In the special case where $\FFF(\omega)$
is always a zero cone for a probability measure on $\Omega$,
this is a protocol for betting on successive outcomes
at odds given by a homogeneous Markov chain on $\Omega$.
However, since there are no probabilities for the initial state $\omega_0$,
our result for this protocol
(part of Theorem \ref{thm:weakly-invariant})
is not directly comparable to the standard zero-one law for Markov chains
mentioned in \S\ref{sec:introduction}.

\medskip

\noindent
\textsc{Protocol 3. Markov trials}

\parshape=2
\IndentI   \WidthI
\IndentII  \WidthII
\noindent
\textbf{Parameters:} set $\Omega$; for each $\omega \in \Omega$,
  a coherent cone $\FFF(\omega)$ of real-valued\\
  functions on $\Omega$

\noindent
\textbf{Protocol:}

\parshape=7
\IndentI   \WidthI
\IndentI   \WidthI
\IndentI   \WidthI
\IndentII  \WidthII
\IndentII  \WidthII
\IndentII  \WidthII
\IndentI   \WidthI
\noindent
Reality announces $\omega_0\in\Omega$.\\
$\K_0 := 1$.\\
FOR $n=1,2,\ldots$:\\
  Skeptic announces $F_n\in\FFF(\omega_{n-1})$.\\
  Reality announces $\omega_n\in\Omega$.\\
  $\K_n := \K_{n-1} + F_n(\omega_n)$.\\
END FOR

\section{Events and upper and lower probabilities}\label{sec:upplowprob}

Upper and lower probabilities can be defined for any of the protocols
used in game-theoretic probability 
(Shafer and Vovk (2001), Takeuchi (2004)). 
We now review the definitions assuming, for simplicity,
that we are using either Protocol 1 or Protocol 2,
where the cone $\FFF_n$ from which Skeptic chooses on each trial is fixed,
independently of how Reality moves earlier in the game.

We call the set $\Omega^{\infty}$ of all infinite sequences of outcomes
the \emph{sample space}.
We write $\omega_1\omega_2\ldots$ for a generic element of $\Omega^\infty$,
and we write $\omega_1\ldots\omega_n$ for a finite sequence of outcomes.

As we will see in this section,
upper and lower probabilities are defined
for any subset $E$ of the sample space $\Omega^{\infty}$.
Accordingly, we call any subset of $\Omega^{\infty}$ an \emph{event}.
This diverges from standard terminology, in which only
elements of a specified $\sigma$-algebra are called events.
Our definitions of \emph{tail event} (\S\ref{sec:kolmogorov}),
\emph{invariant event} (\S\ref{sec:ergodic}),
and \emph{permutable event} (\S\ref{sec:hewitt-savage})
will also make no reference to any $\sigma$-algebra.

A strategy for Skeptic specifies his moves $F_1,F_2,\ldots$
as functions of the preceding moves by Reality:
$F_n$ is a function of $\omega_1 \ldots\omega_{n-1}$.
Once we fix such a strategy, Skeptic's capital process $\K_0,\K_1,\K_2,\ldots$
depends on Reality's moves alone:
$\K_0$ remains the constant $1$,
and $\K_n$ is a function of $\omega_1 \ldots\omega_n$.

We call a strategy for Skeptic \emph{prudent} if its capital process
is everywhere nonnegative---i.e., if
$$
   \K_n(\omega_1\ldots\omega_n) \ge 0
   \quad
   \text{for all $\omega_1\omega_2\ldots\in\Omega^{\infty}$ and all $n$}.
$$
In this case, the strategy risks only the initial unit capital.
A strategy for Skeptic will satisfy
\begin{equation}\label{eq:limit}
  \limsup_{n\to\infty}\K_n(\omega_1\ldots\omega_n)
  \ge 0
  \quad
  \text{for all $\omega_1\omega_2\ldots \in\Omega^{\infty}$}
\end{equation}
if and only if it is prudent.
If the strategy is not prudent---i.e., if
$
  \K_n(\omega_1\ldots\omega_n) < 0
$
for some $\omega_1\omega_2\ldots$ and some $n$,
then Reality can violate~(\ref{eq:limit})
by making Skeptic's subsequent payoffs nonpositive,
which is possible because of the coherence of the cones
from which Skeptic selects his moves.
The $\limsup_{n\to\infty}$ in (\ref{eq:limit})
can be replaced by $\liminf_{n\to\infty}$ or by $\inf_n$.

In our protocols, for any given event $E$ and any $c>0$,
Skeptic has a prudent strategy that satisfies
$$
  \liminf_{n\to\infty}\K_n(\omega_1\ldots\omega_n)
  \ge c
  \quad
  \text{for all $\omega_1 \omega_2 \ldots \in E$}
$$
if and only if he has a prudent strategy that satisfies
\begin{equation*} 
  \sup_{n=1,2,\ldots} \K_n(\omega_1\ldots\omega_n)
  \ge c
  \quad
  \text{for all $\omega_1\omega_2 \ldots \in E$}.
\end{equation*}
This is because Skeptic can stop betting (choose $F_n=0$)
once his capital reaches a given level.

For each event $E$, we set
\begin{multline}\label{eq:success}
  \UpperProb(E)
  :=
  \inf
  \Bigl\{
    \epsilon > 0 \st
    \text{Skeptic has a prudent strategy for which}\\
    \sup_{n=1,2,\ldots} \K_n(\omega_1\ldots\omega_n)
    \ge
    1/\epsilon
    \quad
    \text{for all $\omega_1\omega_2\ldots \in E$}
  \Bigr\},
\end{multline}
and we set
\begin{equation}\label{eq:lower-probability}
  \LowerProb(E)
  :=
  1
  -
  \UpperProb
  \left(
    E^c
  \right).
\end{equation}
We call $\UpperProb(E)$ the \emph{upper probability} of $E$, and we call
$\LowerProb(E)$ the \emph{lower probability} of $E$.
Upper and lower probability generalize the standard notion of probability.
Let us consider, for simplicity, Protocol 1.
If $\FFF$ is a zero cone for a probability measure $P$
on a $\sigma$-algebra for $\Omega$,
then $\UpperProb(E)=\LowerProb(E)=P^{\infty}(E)$
for all measurable subsets of $\Omega$ 
(Shafer and Vovk (2001), \S 8.2),
and our theorems below reduce to the standard zero-one laws
in measure-theoretic probability.

\begin{lemma}\label{lem:coherence}
  $0 \le \LowerProb(E)\le \UpperProb(E) \le 1$
  for every event $E$.
\end{lemma}

\begin{proof}
  The relation $\UpperProb(E) \le 1$ follows from the fact that Skeptic
  can choose $F_n=0$ for all $n$, and $0 \le \LowerProb(E)$ then follows
  by the definition~(\ref{eq:lower-probability}).

  Suppose $\LowerProb(E)>\UpperProb(E)$,
  i.e.,
  $\UpperProb(E)+\UpperProb(E^c)<1$.
  Then there exist $\epsilon_1>0$, $\epsilon_2>0$,
  and prudent strategies $S_1$ and $S_2$ for Skeptic
  such that $\epsilon_1+\epsilon_2<1$,
  $S_1$ guarantees $\sup_n\K_n\ge\III_E/\epsilon_1$,
  and $S_2$ guarantees $\sup_n\K_n\ge\III_{E^c}/\epsilon_2$,
  where $\III_E$ is the indicator function of $E$, i.e.,
  $\III_E:\Omega^{\infty}\to\bbbr$ takes the value $1$ on $E$
  and the value $0$ outside $E$.
  Then $(\epsilon_1S_1+\epsilon_2S_2)/(\epsilon_1+\epsilon_2)$
  guarantees
  \[
    \sup_n\K_n\ge\ \frac{\III_{E\cup E^c}}{\epsilon_1+\epsilon_2}
    =
    \frac{1}{\epsilon_1+\epsilon_2}
    >
    1.
  \]
  But this is impossible
  since, by coherence,
  Reality can choose $\omega_1 \omega_2 \ldots$
  so that $1=\K_0\ge\K_1\ge\K_2\ge\cdots$.
\end{proof}

We can also consider capital processes determined by different strategies for Skeptic
when his initial capital $\K_0$ is not necessarily equal to $1$.
We call any such capital process a \emph{martingale}.
The martingales form a cone.
We can rephrase our definition of upper probability, (\ref{eq:success}),
by saying that $\UpperProb(E)$ is the infimum of all values of $\epsilon$
such that there exists a nonnegative martingale starting at $\epsilon$
and reaching at least $1$ on every sequence $\omega_1\omega_2\ldots$ in $E$.

The preceding definitions are easily adapted to Protocol~3;
we simply recognize that Skeptic's strategies and martingales
will also depend on $\omega_0$
as well as on $\omega_1\omega_2\ldots$\,.
These definitions also apply, with similarly minor modifications,
to other protocols used in game-theoretic probability.

The following terminology spells out the intuitive meaning
of extreme values for upper and lower probabilities:
\begin{itemize}
  \item
     When $\LowerProb(E)=0$, we say $E$ is \emph{unsupported}.
  \item
     When $\LowerProb(E)=1$, we say $E$ is \emph{certain}.
  \item
     When $\UpperProb(E)=0$, we say $E$ is \emph{impossible}.
  \item
     When $\UpperProb(E)=1$, we say $E$ is \emph{fully plausible}.
  \item
     When $E$ is unsupported and fully plausible,
     we say it is \emph{fully uncertain}.
\end{itemize}
When an event is certain, it is also fully plausible.
When it is impossible, it is also unsupported.
An event being unsupported is equivalent
to its complement being fully plausible.
An event being certain is equivalent to its complement being impossible.
An event is fully uncertain if and only if its complement is fully uncertain.

As we remarked in \S\ref{sec:introduction}, both $\LowerProb(E)$ and $\UpperProb(E)$
will coincide with $E$'s probability
when there are enough prices to determine a probability measure
for a $\sigma$-algebra containing $E$;
see Shafer and Vovk (2001), 
\S 8.2.
In this case, $E$ cannot be fully uncertain.

\section{Two types of zero-one law}\label{sec:twotypes}

A measure-theoretic zero-one law says that an event $E$
satisfying specified conditions is either impossible or certain:
either $\Prob(E)=0$ or $\Prob(E)=1$.
In the general game-theoretic case,
where we have only upper and lower probabilities,
we get one of the following weaker statements:
\begin{enumerate}
  \item
    $E$ is either fully plausible or unsupported (or both---i.e., fully uncertain).
  \item
    $E$ is certain, impossible, or fully uncertain.
\end{enumerate}
Condition~2 is stronger than Condition~1,
because certain implies fully plausible,
and impossible implies unsupported.

When the prices determine a probability measure
on a $\sigma$-algebra containing $E$,
we have $\UpperProb(E)=\LowerProb(E)=\Prob(E)$,
and both conditions then imply that $\Prob(E)=1$ or $\Prob(E)=0$.

Our game-theoretic versions of Kolmogorov's zero-one law
and ergodicity will assert Condition~2,
but our game-theoretic version of Hewitt and Savage's zero-one law
(proven only in a special case)
will assert only Condition~1.

\section{Kolmogorov's zero-one law}
\label{sec:kolmogorov}

Our game-theoretic version of Kolmogorov's zero-one law is a theorem about Protocol~2. The relevant events are tail events.

An event $E$ 
is  called a \emph{tail event}
if any sequence in $\Omega^{\infty}$ that agrees from some point onwards
with a sequence in $E$ is also in $E$---i.e.,
if $\omega_1\omega_2\ldots$ and $\omega_1^{\prime} \omega_2^{\prime}\ldots$
are either both in $E$ or both not in $E$ whenever
$\omega_n=\omega_n^{\prime}$ except for a finite number of $n$.
It follows immediately from this definition that
$E$ is a tail event if and only if its complement $E^c$ is a tail event.

\begin{theorem}\label{thm:kolmogorov}
Suppose that $E$ is a tail event in Protocol~2 (independently priced trials).
Then $E$ is certain, impossible, or fully uncertain.
\end{theorem}

Skeptic chooses from $\FFF_n$ on trial $n$ in Protocol~2.
Our proof of Theorem~\ref{thm:kolmogorov} will use
the fact that we get another instantiation of Protocol~2
if we start on trial $n+1$
for some $n \ge 1$---i.e.,
if Skeptic chooses from $\FFF_{n+1}$ on the first trial,
from $\FFF_{n+2}$ on the second trial, and so on.
We call this the \emph{shifted protocol}, as opposed to the \emph{original protocol}.
The two protocols, the original one and the shifted one, have the same events;
in both, any subset of $\Omega^{\infty}$ is an event.

Let $\UpperProb^{-n}$ denote upper probability in the shifted protocol.
(In particular, $\UpperProb^{-0} = \UpperProb$.)
Given a strategy $S$ for Skeptic in the shifted protocol,
write $S^{+n} $ for the strategy in the original protocol
that sets Skeptic's first $n$ moves equal to $0$ and then plays $S$.

We write $\theta$ for the shift operator,
which deletes the first element from a sequence in $\Omega^{\infty}$:
\[
  \theta: \omega_1 \omega_2 \omega_3 \ldots \mapsto \omega_2 \omega_3 \ldots\,.
\]
We write $E^{-n}$ for $\theta^n E$:
\begin{equation*} 
  E^{-n}
  :=
  \{\  \omega_{n+1} \omega_{n+2} \ldots\  
  \st \ \omega_1 \omega_2\ldots \in E \}.
\end{equation*}

The next two lemmas relate upper probabilities
in the original and shifted protocols.
\begin{lemma}
  Let $E$ be any event in Protocol~2.  Then 
 $\UpperProb^{-n}(E^{-n})$ is non-decreasing in $n$, i.e.,
  $\UpperProb(E)\le\UpperProb^{-1}(E^{-1})\le\UpperProb^{-2}(E^{-2})\le\cdots$.
\end{lemma}
\begin{proof}
Suppose $c > 0$,
and suppose $S$ is a prudent strategy in the shifted protocol that achieves
$$
  \sup_{k=1,2,\ldots} \K_k(\omega_1\ldots\omega_k)
  \ge c
  \quad
  \text{for all $\omega_1\omega_2\ldots \in E^{-n}$}.
$$
Then $S^{+n}$ is evidently also prudent and achieves
$$
  \sup_{k=1,2,\ldots} \K_k(\omega_1\ldots\omega_k)
  \ge c
  \quad
  \text{for all $\omega_1\omega_2\ldots \in E$}
$$
in the original protocol.
Let ${\cal S}^{+n}$ denote
the set of strategies in the original protocol of the form  $S^{+n}$.
It follows that the infimum in (\ref{eq:success}) over ${\cal S}^{+n}$
coincides with $\UpperProb^{-n}(E^{-n})$.
Since  the set ${\cal S}^{+n}$ is non-increasing in $n$,
we are taking the infimum over a smaller set of strategies in (\ref{eq:success})
when $n$ is larger.  So
$\Prob^{-n}(E^{-n})$ is non-decreasing in $n$.
\end{proof}

\begin{lemma}\label{lem:tail-event}
  Suppose that $E$ is a tail event in Protocol~2. Then
  $
    \UpperProb(E) = \UpperProb^{-n}(E^{-n})
  $
  for all $n$.
\end{lemma}
\begin{proof}
Suppose $c > 0$,
and suppose $S$ is a prudent strategy in the original protocol
that achieves
\begin{equation}\label{eq:finally}
  \sup_{k=1,2,\ldots} \K_k(\omega_1 \ldots\omega_k)
  \ge c
  \quad
  \text{for all $\omega_1 \omega_2\ldots \in E$}.
\end{equation}
Because $E$ is a tail event,
Reality can choose any sequence from $\Omega^n$
as her first $n$ moves without affecting whether $E$ happens.
By coherence, she can choose these $n$ moves
so that $S$ makes no money for Skeptic on the $n$ trials.
Condition~(\ref{eq:finally}) tells us
that the moves specified by $S$ on the $(n+1)$th and later trials
guarantee that Skeptic can make up this loss
and still get at least $c$
when $\omega_1\omega_2\ldots\in E$.
In the shifted protocol, Skeptic has capital $1$ at the beginning
rather than the same or smaller capital resulting from the losses.
So these moves still define a prudent strategy and guarantee that
$$
  \sup_{k=1,2,\ldots} \K_k(\omega_1\ldots\omega_k)
  \ge c
  \quad
  \text{for all $\omega_1 \omega_2\ldots \in E^{-n}$}
$$
in the shifted protocol.
This implies
$\UpperProb^{-n}(E^{-n})\le \UpperProb(E)$
and together with the previous lemma
we obtain $\UpperProb^{-n}(E^{-n})=\UpperProb(E)$.
\end{proof}

\begin{proof}[Proof of Theorem~\ref{thm:kolmogorov}]
  We will show that if a tail event is not fully plausible,
  then it is impossible.
  This suffices to prove the theorem,
  because if $E$ is a tail event, then $E^c$ is also a tail event.
  If $E$ is not fully uncertain, then either $E$
  or $E^c$ is not fully plausible, and if one of them is impossible,
  then $E$ is either impossible or certain.

  Suppose, then, that $E$ is not fully plausible: $\UpperProb(E)<1$.
  Assume that Reality chooses a path in $E$.
  Choose $\epsilon$ such that $\UpperProb(E)<\epsilon<1$.
  Then there is a prudent strategy $S$ for Skeptic that guarantees
  he will multiply his initial capital of $1$ by $1/\epsilon$
  in a finite number of trials.
  Skeptic plays this strategy until his capital reaches at least $1/\epsilon$.
  Then he starts over, playing $(1/\epsilon)S$ in the shifted game starting at that point.
  This eventually again multiplies his capital by another factor of $1/\epsilon$ or more.
  Continuing in this way,
  he can make his capital arbitrarily large while playing prudently.
  This demonstrates that $\UpperProb(E)=0$---i.e., that $E$ is impossible.
\end{proof}

We have just shown that when $E$ is a tail event
with $\UpperProb(E)<1$,
\begin{equation}\label{eq:strongly}
  \text{Skeptic has a prudent strategy guaranteeing }
    \lim_{n\to\infty}\K_n=\infty \text{ on } E.
\end{equation}
This implies $\UpperProb(E)=0$, but in general it may be stronger.
So we have proven a bit more than the theorem asserts.
We have proven that a tail event is either strongly certain,
strongly impossible, or fully uncertain,
where an event is said to be \emph{strongly impossible}
if~(\ref{eq:strongly}) holds
and \emph{strongly certain} if its complement is strongly impossible.

Here is an example of a fully uncertain tail event.
Consider Protocol~2 where $\Omega$ is a linear space,
$\Omega\ne\{0\}$, and $\FFF_n \equiv \FFF$
is the set of linear functions on $\Omega$. If Reality chooses the origin
$0\in \Omega$, then $F(0)=0$ for all $F \in \FFF$.
On the other hand
if Skeptic chooses $F\neq 0$,
then there exists $\omega\in\Omega$ such that $F(\omega) < 0$.
Therefore the protocol is coherent.
Let $E$ be the event that
Reality chooses $\omega_n=0$ except for a finite number of $n$.
Clearly $E$ is a tail event.
Since $000\ldots\in E$, $\UpperProb(E)=1$.
Now consider $E^c$,
which is the event that $\omega_n \neq 0$ infinitely often.
For each choice $F\in \FFF$,
Reality can choose $\omega \neq 0$ such that $F(\omega)\le 0$.
So $\UpperProb(E^c)=1$.

Suppose Reality is determined to prevent Skeptic's capital
from becoming arbitrarily large.
If $E$ is a fully uncertain tail event,
Skeptic has no control at all over Reality choosing an element from the event $E$.
On the other hand, there are many examples of tail events
in Shafer and Vovk (2001) 
which are certain (i.e., Skeptic can force these events),
such as the strong law of large numbers or the law of the iterated logarithm.

It is interesting to compare our martingale proof
with Kolmogorov's measure-theoretic proof,
given in an appendix to \emph{Grundbegriffe}
and reproduced in many textbooks.
Kolmogorov shows that a tail event $E$ is independent of itself,
so that $\Prob(E)=\Prob(E)^2$ and therefore
$\Prob(E)=0 \text{ or }\Prob(E)=1$.
Our martingale proof paints a little larger picture.
Having Skeptic start over just once
after multiplying his capital by $1/\epsilon$
suffices to show that $\UpperProb(E)^2\ge\UpperProb(E)$,
and this implies $\UpperProb(E) =0$ or $\UpperProb(E)=1$.
But by having Skeptic start over again and again,
we find that he can become infinitely rich
if an event $E$ with $\UpperProb(E)<1$ happens---i.e.,
that such $E$ is strongly impossible.

A proof based on the measure-theoretic martingale convergence theorem
is also known
(\S 2.2 of Chow et al.\ (1971), \S 14.3 of Williams (1991)),
but it is much less constructive
and restricted to measurable events.
A corresponding constructive game-theoretic proof  is desirable.
A referee pointed out a result for tail events by B\'artfai and R\'ev\'esz 
(B\'artfai and R\'ev\'esz (1967), Theorems 1 and 2)
which can be interpreted as an approximate zero-one law
that holds under weak long-range dependence.
It would be interesting to establish a similar result
in the game-theoretic framework by appropriately modifying our protocols.

\section{Ergodicity}
\label{sec:ergodic}

Now consider Protocol~1, where Skeptic always chooses from the same cone,
and Protocol~3, where the cone may depend on Reality's previous move.
The relevant events are invariant events.

We call an event $E$ 
\emph{weakly invariant}
if $\theta E = E^{-1}\subseteq E$.

If $E$ is weakly invariant, then by induction $E^{-n}$ is non-increasing in $n$.
In accordance with standard terminology 
(e.g., Shiryaev (1996), \S V.2),
we call an event $E$ \emph{invariant} if $E=\theta^{-1}E$.

\begin{lemma}\label{lem:weakly-invariant}
  $E$ is invariant if and only if both $E$ and $E^c$ are weakly invariant.
\end{lemma}

\begin{proof}
If $E$ is invariant, then $E^c$ is also invariant,
because the inverse map commutes with complementation.
Hence in this case both $E$ and $E^c$ are weakly invariant.

Conversely suppose that $\theta E \subseteq E$ and $\theta E^c \subseteq E^c$.
The first inclusion is equivalent to $E\subseteq\theta^{-1}E$
and the second is equivalent to $E^c\subseteq\theta^{-1}E^c$.
Since the right-hand sides of the last two inclusions
are disjoint,
these inclusions are in fact equalities.
\end{proof}

\begin{theorem}\label{thm:weakly-invariant}
  Suppose that $E$ is a weakly invariant event in Protocol~1 (identically priced trials)
  or Protocol~3 (Markov trials).  Then $E$ is either impossible or fully plausible.
\end{theorem}

\begin{proof}
  It suffices, as in the proof of Theorem~\ref{thm:kolmogorov},
  to show that $\UpperProb(E)<1$ implies $\UpperProb(E)=0$.
  The current proof is, however,
  significantly simpler than that of Theorem~\ref{thm:kolmogorov}:
  no analogue of Lemma \ref{lem:tail-event} is needed.

  Let $E$ be a weakly invariant event in Protocol~1
  such that $\UpperProb(E)< \epsilon <1$.
  Assume that Reality chooses a path in $E$.
  Skeptic plays a prudent strategy until his capital increases
  by a factor of $1/\epsilon$ or more at some round $n_1$.
  By the assumption of weak invariance, $E^{-n_1} \subseteq E$.
  This implies that when he starts over,
  there will be another time $n_2>n_1$
  such that he multiplies his capital again by $1/\epsilon$ or more.
  Continuing in this way,
  he can make his capital arbitrarily large while playing prudently.

  In the case of Protocol~3, Skeptic's strategies and martingales
  depend on $\omega_0$ as well as on $\omega_1\omega_2\ldots$,
  and the definition of $\UpperProb(E)$, (\ref{eq:success}), becomes
  \begin{multline*}
    \UpperProb(E)
    :=
    \inf
    \Bigl\{
      \epsilon > 0
      \st
      \text{Skeptic has a prudent strategy for which}\\
      \sup_{n=1,2,\ldots} \K_n(\omega_0 \omega_1 \ldots\omega_n)
      \ge
      1/\epsilon
      \quad
      \text{for all $\omega_0\omega_1\omega_2\ldots \in E$}
    \Bigr\}.
  \end{multline*}
  Similarly,
  for each $\omega\in\Omega$ we define
  \begin{multline*}
    \UpperProb_{\omega}(E)
    :=
    \inf
    \Bigl\{
      \epsilon > 0
      \st
      \text{Skeptic has a prudent strategy for which}\\
      \sup_{n=1,2,\ldots} \K_n(\omega\omega_1\ldots\omega_n)
      \ge
      1/\epsilon
      \quad
      \text{for all $\omega_1\omega_2\ldots\in E$}
    \Bigr\}.
  \end{multline*}
  We have
  $
    \UpperProb(E) = \sup_{\omega\in\Omega} \UpperProb_{\omega}(E)
  $.
  Suppose $\UpperProb(E)<\epsilon<1$.
  For each $\omega\in\Omega$,
  consider the modified protocol
  in which Reality's first move is $\omega$,
  and choose a prudent strategy $S_{\omega}$ for Skeptic in this protocol
  that eventually multiplies the initial unit capital by $1/\epsilon$ or more.

  Now we can proceed as before:
  Skeptic first plays $S_{\omega_0}$ until the first $n$ for which $\K_n \ge 1/\epsilon$,
  then plays a scaled up version of $S_{\omega_{n}}$
  until his capital is again multiplied by $1/\epsilon$ or more, etc.
\end{proof}

In view of Lemma \ref{lem:weakly-invariant}
we obtain the following corollary to Theorem \ref{thm:weakly-invariant}.

\begin{corollary}\label{thm:bernoulli}
  Suppose that $E$ is an invariant event in Protocol~1 (identically priced trials)
  or Protocol~3 (Markov trials).
  Then $E$ is certain, impossible, or fully uncertain.
\end{corollary}

\section{Hewitt and Savage's zero-one law}
\label{sec:hewitt-savage}

Let us again consider Protocol~1, where Skeptic always chooses from the same cone.

Let us call an event $E$ 
\emph{permutable}
if for any sequence in $E$ and any $n>1$,
any sequence obtained by permuting the first $n$ terms of the sequence
is also in $E$.
Let us call $E$ \emph{singly generated}
if it is equal to the set of sequences obtained by taking a single sequence
and permuting finite initial subsequences in all possible ways.

We may conjecture that any permutable event in Protocol~1
is either fully plausible or unsupported.
But we can prove this result
only in the case where the permutable event is singly generated.

\begin{proposition}\label{prop:hewitt-savage-one-sequence}
  Suppose that $E$ is a singly generated permutable event in Protocol~1
  (identically priced trials). 
  Then $E$ is either fully plausible or unsupported.
\end{proposition}
\begin{proof}[Proof of Proposition~\ref{prop:hewitt-savage-one-sequence}]
  Suppose $E$ is neither fully plausible nor unsupported.
  Then we may choose $\epsilon<1$ such that
  $\UpperProb(E)<\epsilon$ and $\UpperProb(E^c)<\epsilon$.
  Choose prudent strategies $S$ and $S'$ for Skeptic
  that multiply his capital by at least $1/\epsilon$
  on $E$ and $E^c$, respectively.

  As usual, we derive a contradiction by showing
  the existence of a prudent strategy for Skeptic
  that makes his capital tend to infinity if $E$ happens
  (i.e., showing that $E$ is strongly impossible and so unsupported).
  Let $\omega_1\omega_2\ldots$ be the sequence of outcomes
  actually chosen by Reality.
  The strategy begins by playing $S$ until the capital exceeds $1/\epsilon$.
  This happens on some trial $n_1$ if $\omega_1\omega_2\ldots \in E$.
  At this point,
  we ask whether the remaining sequence
  $\omega_{n_1+1}\omega_{n_1+2}\ldots$
  is in $E$ or $E^c$.
  In the first case, the strategy plays a scaled up version of $S$;
  in the second it plays a scaled up version of $S'$.
  Again, the capital will eventually be multiplied by $1/\epsilon$
  on some trial $n_2$.
  Etc.

  In conclusion,
  let us check that the initial sequence $\omega_1\ldots\omega_n$
  (where $n\in\{n_1,n_2,\ldots\}$)
  indeed determines whether the remaining sequence
  $\omega_{n+1}\omega_{n+2}\ldots$
  is in $E$ or $E^c$.
  Suppose there are two possible continuations
  \begin{equation}\label{eq:continuations}
    \omega'_{n+1}\omega'_{n+2}\ldots\in E
    \text{\quad and\quad}
    \omega''_{n+1}\omega''_{n+2}\ldots\notin E
  \end{equation}
  such that both
  $\omega_1\ldots\omega_n\omega'_{n+1}\omega'_{n+2}\ldots$
  and
  $\omega_1\ldots\omega_n\omega''_{n+1}\omega''_{n+2}\ldots$
  belong to $E$.
  Since $E$ is singly generated,
  for some $N>n$ it is true that:
  (a) the sequences
  $\omega_1\ldots\omega_n\omega'_{n+1}\ldots\omega'_{N}$
  and
  $\omega_1\ldots\omega_n\omega''_{n+1}\ldots\omega''_{N}$
  are permutations of each other;
  (b) $\omega'_i=\omega''_i$ for $i>N$.
  Condition (a) means that the corresponding multisets,
  $\lbag\omega_1,\ldots,\omega_n,\omega'_{n+1},\ldots,\omega'_{N}\rbag$
  and
  $\lbag\omega_1,\ldots,\omega_n,\omega''_{n+1},\ldots,\omega''_{N}\rbag$,
  coincide.
  (We write $\lbag\ldots\rbag$ rather than $\{\ldots\}$
  to emphasize that repetitions are allowed.)
  Therefore,
  the multisets
  $\lbag\omega'_{n+1},\ldots,\omega'_{N}\rbag$
  and
  $\lbag\omega''_{n+1},\ldots,\omega''_{N}\rbag$
  also coincide.
  In other words,
  the sequences
  $\omega'_{n+1}\ldots\omega'_{N}$
  and
  $\omega''_{n+1}\ldots\omega''_{N}$
  are permutations of each other.
  In combination with condition (b),
  this contradicts our assumption (\ref{eq:continuations}).
\end{proof}

Here is an example showing
that the conclusion of Proposition~\ref{prop:hewitt-savage-one-sequence}
cannot be strengthened to say that $E$ is impossible or fully plausible;
in particular, to say that $E$ is certain, impossible, or fully uncertain.
Let $\Omega:=\{-1,0,1\}$,
let $\FFF$ consist of functions $F$ taking values $-t,0,t$ at $-1,0,1$,
respectively,
for some $t\ge0$,
and let $E$ be the set of all sequences in $\Omega^{\infty}$
that do not contain $-1$ and contain precisely one $1$.
Then $E$ is permutable and singly generated;
it is generated by the sequence $100\ldots$\,.
Reality can keep Skeptic from making any money
by choosing the sequence $000\ldots$,
and since this sequence is in $E^c$,
this implies that $\UpperProb(E^c) = 1$ and $\LowerProb(E)=0$.
This is consistent with the conclusion of the proposition:
$E$ is unsupported.
But the prudent strategy that does best for Skeptic on $E$
is one that chooses $t=1$ at the first trial
and continues with this choice so long as Reality plays $0$;
this doubles Skeptic's money on $E$ but no more,
and so $\UpperProb(E)=1/2$.
Thus $E$ is neither impossible nor fully plausible.

It would be interesting
to extend Proposition \ref{prop:hewitt-savage-one-sequence}
to all permutable events
or construct a permutable event $E$ for which
$
  0
  <
  \LowerProb(E)
  \le
  \UpperProb(E)
  <
  1.
$
Note that Proposition \ref{prop:hewitt-savage-one-sequence}
can be trivially extended to the case of permutable events
generated by two or more sequences that are distinguishable by any initial segment,
e.g., generated by a sequence consisting of rational numbers
and a sequence consisting of irrational numbers.

\bibliographystyle{plain}

\end{document}